# Optimal Plant Layout Design for Process-focused Systems

*Corresponding Author: Dr. M. Khoshnevisan*

*School of Accounting & Finance*

*Griffith University, Australia*

*E-mail: m.khoshnevisan@mailbox.gu.edu.au*

*Tel (0061)7-555-28763*

*Sukanto Bhattacharya*

*School of Information Technology*

*Bond University, Australia*

*Dr. Florentin Smarandache*

*University of New Mexico - Gallup, USA*

## Abstract

In this paper we have proposed a semi-heuristic optimization algorithm for designing optimal plant layouts in process-focused manufacturing/service facilities. Our proposed algorithm marries the well-known CRAFT (Computerized Relative Allocation of Facilities Technique) with the Hungarian assignment algorithm. Being a semi-heuristic search, our algorithm is likely to be more efficient in terms of computer CPU engagement time as it tends to converge on the global optimum faster than the traditional CRAFT algorithm - a pure heuristic. We also present a numerical illustration of our algorithm.

# Key Words

Facilities layout planning, load matrix, CRAFT, Hungarian assignment algorithm

**Introduction**

The fundamental integration phase in the design of productive systems is the layout of production facilities. A working definition of layout may be given as the arrangement of machinery and flow of materials from one facility to another, which minimizes material-handling costs while considering any physical restrictions on such arrangement.

Usually this layout design is either on considerations of machine-time cost and product availability; thereby making the production system *product-focused*; or on considerations of quality and flexibility; thereby making the production system *process-focused*.

It is natural that while *product-focused* systems are better off with a 'line layout' dictated by available technologies and prevailing job designs, *process-focused* systems, which are more concerned with job organization, opt for a 'functional layout'. Of course, in reality the actual facility layout often lies somewhere in between a pure line layout and a pure functional layout format; governed by the specific demands of a particular production plant. Since our present paper concerns only functional layout design for process-focused systems, this is the only layout design we will discuss here.

The main goal to keep in mind is to minimize material handling costs - therefore the departments that incur the most interdepartmental movement should be located closest to one another. The main type of design layouts is *Block diagramming*, which refers to the movement of materials in existing or proposed facility. This information is usually provided with a *from/to* chart or *load summary* chart, which gives the average number of units loads moved between departments. A load-unit can be a single unit, a pallet of material, a bin of material, or a crate of material. The next step is to design the layout by calculating the composite movements between departments and rank them from most movement to least movement. Composite movement refers to the back-and-forth movement between each pair of departments. Finally, trial layouts are place on a grid that graphically represents the relative distances between departments. This grid then becomes the objective of optimization when determining the optimal plant layout.

We give a visual representation of the basic operational considerations in a *process-focused* system schematically as follows:

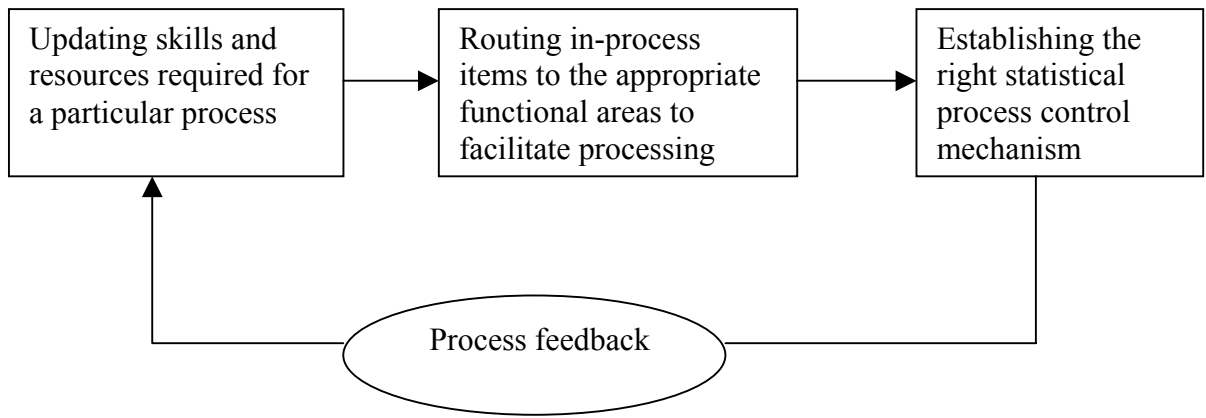

In designing the optimal functional layout, the fundamental question to be addressed is that of 'relative location of facilities'. The locations will depend on the need for one pair of facilities to be adjacent (or physically close) to each other relative to the need for all other pairs of facilities to be similarly adjacent (or physically close) to each other. Locations must be allocated based on the relative gains and losses for the alternatives and seek to minimize some indicative measure of the cost of having non-adjacent locations of facilities. Constraints of space prevents us from going into the details of the several criteria used to determine the gains or losses from the relative location of facilities and the available sequence analysis techniques for addressing the question; for which we refer the interested reader to any standard handbook of production/operations management.

**Computerized Relative Allocation of Facilities Technique (CRAFT)**

CRAFT (Buffa et al., 1964) is a computerized heuristic algorithm that takes in *load matrix* of interdepartmental flow and transaction costs with a representation of a block layout as the inputs. The block layout could either be an existing layout or; for a new facility, any arbitrary initial layout. The algorithm then computes the departmental locations and returns an estimate of the total interaction costs for the initial layout. The governing algorithm is designed to compute the impact on a cost measure for two-way or three-way swapping in the location of the facilities. For each swap, the various interaction costs are computed afresh and the load matrix and the change in cost (increase or decrease) is noted and stored in the RAM. The algorithm proceeds this way through all possible combinations of swaps accommodated by the software. The basic procedure is repeated a number of times resulting in a more efficient block layout every time till such time when no further cost reduction is possible. The final block layout is then printed out to serve as the basis for a detailed layout template of the facilities at a later stage. Since its formulation, more powerful versions of CRAFT have been developed but these too follow the same, basic heuristic routine and therefore tend to be highly CPU-intensive.

The basic computational disadvantage of a CRAFT-type technique is that one always has got to start with an arbitrary initial solution. This means that there is no mathematical certainty of attaining the desired optimal solution after a given number of iterations. If the starting solution is quite close to the optimal solution by chance, then the final solution is attained only after a few iterations. However, as there is no guarantee that the starting solution will be close to the global optimum, the expected number of iterations required to arrive at the final solution tend to be quite large thereby straining computing resources.

In our present paper we propose and illustrate the Modified Assignment (MASS) algorithm as an extension to the traditional CRAFT, to enable faster convergence to the optimal solution. This we propose to do by marrying CRAFT technique with the Hungarian assignment algorithm. As our proposed algorithm is semi-heuristic, it is likely to be less CPU-intensive than any traditional, purely heuristic CRAFT-type algorithm.

**The Hungarian assignment algorithm**

A general assignment problem may be framed as a special case of the *balanced transportation problem* with availability and demand constraints summing up to unity. Mathematically, it has the following *general linear programming* form:

**Minimize $\Sigma\Sigma$ $C_{ij}X_{ij}$**

**Subject to $\Sigma X_{ij} = 1$, for each i, j = 1, 2 …n .**

In words, the problem may be stated as assigning each of n individuals to n jobs so that exactly one individual is assigned to each job in such a way as to minimize the total cost.

To ensure satisfaction of the basic requirements of the assignment problem, the basic feasible solutions of the corresponding balanced transportation problem must be integer valued. However, any such basic feasible solution will contain (2n – 1) variables out of which (n – 1) variables will be zero thereby introducing a high level of degeneracy in the solution making the usual solution technique of a transportation problem very inefficient. This has resulted in mathematicians devising an alternative, more efficient algorithm for solving this class of problems, which has come to be commonly known as the *Hungarian assignment algorithm*. Basically, this algorithm draws from a simple theorem in linear algebra which says that if a constant number is added to any row and/or column of the cost matrix of an assignment-type problem, then the resulting assignment-type problem has exactly the same set of optimal solutions as the original problem and vice versa.

**Proof:**

Let $A_i$ and $B_j$ **(i, j = 1, 2 … n)** be added to the ith row and/or jth column respectively of the cost matrix. Then the revised cost elements are $C_{ij}^* = C_{ij} + A_i + B_j$. The revised cost of assignment is $\Sigma\Sigma C_{ij}^* X_{ij} = \Sigma\Sigma (C_{ij} + A_i + B_j) X_{ij} = \Sigma\Sigma C_{ij} X_{ij} + \Sigma A_i \Sigma X_{ij} + \Sigma B_j \Sigma X_{ij}$. But by the imposed assignment constraint $\Sigma X_{ij} = 1$ **(for i, j = 1, 2 … n)**, we have the revised cost as $\Sigma\Sigma C_{ij} X_{ij} + \Sigma A_i + \Sigma B_j$ i.e. the cost differs from the original by a constant. As the revised costs differ from the originals by a constant, which is independent of the decision variables, an optimal solution to one is also optimal solution to the other and vice versa. This theorem can be used in two different ways to solve the assignment problem. First, if in an assignment problem, some cost elements are negative, the problem may be converted into an equivalent assignment problem by adding a positive constant to each of the entries in the cost matrix so that they all become non-negative. Next, the important thing to look for is a feasible solution that has zero assignment cost after adding suitable constants to the rows and columns. Since it has been assumed that all entries are now non-negative, this assignment must be the globally optimal one.

Given a zero assignment, a straight line is drawn through it (a horizontal line in case of a row and a vertical line in case of a column), which prevents any other assignment in that particular row/column. The governing algorithm then seeks to find the minimum number of such straight lines, which would cover all the zero entries to avoid any redundancy. Let us say that *k* such lines are required to cover all the zeroes. Then the *necessary condition* for optimality is that number of zeroes assigned is equal to *k* and the *sufficient condition* for optimality is that *k* is equal to *n* for an *n x n* cost matrix.

**The MASS (Modified Assignment) algorithm**

The basic idea of our proposed algorithm is to develop a systematic scheme to arrive at the initial input block layout to be fed into the CRAFT program so that the program does not have to start off from any initial (and possibly inefficient) solution. Thus, by subjecting the problem of finding an initial block layout to a mathematical scheme, we in effect reduce the purely heuristic algorithm of CRAFT to a semi-heuristic one. Our proposed MASS algorithm follows the following *sequential* steps:

**Step 1:** We formulate the load matrix such that each entry $l_{ij}$ represents the load carried from facility i to facility j

**Step 2:** We insert $l_{ij}$ = M, where M is a large positive number, into all the vacant cells of the load matrix signifying that no inter-facility load transportation is required or possible between the $i^{th}$ and $j^{th}$ vacant cells

**Step 3:** We solve the problem on the lines of a standard assignment problem using the Hungarian assignment algorithm treating the load matrix as the cost matrix

**Step 4:** We draft the initial block layout trying to keep the inter-facility distance $d_{ij}^*$ between the $i^{th}$ and $j^{th}$ *assigned* facilities to the minimum possible magnitude, subject to the available floor area and architectural design of the shop floor

**Step 5:** We proceed using the CRAFT program to arrive at the optimal layout by iteratively improving upon the starting solution provided by the Hungarian assignment algorithm till the overall load function **L = ΣΣ $l_{ij}d_{ij}^*$** subject to any particular bounds imposed on the problem

The Hungarian assignment algorithm will ensure that the initial block layout is at least very close to the global optimum if not globally optimal itself. Therefore the subsequent CRAFT procedure will converge on the global optimum much faster starting from this near-optimal initial input block layout and will be much less CPU-intensive that any traditional CRAFT-type algorithm. Thus MASS is not a stand-alone optimization tool but rather a rider on the traditional CRAFT that tries to ensure faster convergence to the optimal block layout for process-focused systems, by making the search semi-heuristic.

We provide a numerical illustration of the MASS algorithm in the Appendix by designing the optimal block layout of a small, single-storied, process-focused manufacturing plant with six different facilities and a rectangular shop floor design. The model can however be extended to cover bigger plants with more number of facilities. Also the MASS approach we have advocated here can even be extended to deal with the multi-floor version of CRAFT (Johnson, 1982) by constructing a separate assignment table for each floor subject to any predecessor-successor relationship among the facilities.

*******

# Appendix: Numerical illustration of MASS

We consider a small, single-storied process-focused manufacturing plant with a rectangular shop floor plan having six different facilities. We mark these facilities as $F_I$, $F_{II}$, $F_{III}$, $F_{IV}$, $F_V$ and $F_{VI}$. The architectural design requires that there be an aisle of at least 2meters width between two adjacent facilities and the total floor area of the plant is 64meters x 22meters. Based on the different types of jobs processed, the loads to be transported between the different facilities are supplied in the following load matrix:

|         | $F_I$ | $F_{II}$ | $F_{III}$ | $F_{IV}$ | $F_V$ | $F_{VI}$ |
|---------|-------|----------|-----------|----------|-------|----------|
| $F_I$   | –     | 20       | –         | –        | –     | 25       |
| $F_{II}$ | 10   | –        | 15        | –        | –     | –        |
| $F_{III}$ | –  | –        | –         | 30       | –     | –        |
| $F_{IV}$ | –   | –        | 50        | –        | –     | 40       |
| $F_V$   | –     | –        | –         | –        | –     | 10       |
| $F_{VI}$ | –   | –        | –         | –        | 15    | –        |

We put in a very large positive value M in each of the vacant cells of the load matrix to signify that no inter-facility transfer of load is required or is permissible for these cells:

|       | F$_I$ | F$_{II}$ | F$_{III}$ | F$_{IV}$ | F$_V$ | F$_{VI}$ |
|-------|-------|----------|-----------|----------|-------|----------|
| F$_I$    | M  | 20 | M  | M  | M  | 25 |
| F$_{II}$  | 10 | M  | 15 | M  | M  | M  |
| F$_{III}$ | M  | M  | M  | 30 | M  | M  |
| F$_{IV}$  | M  | M  | 50 | M  | M  | 40 |
| F$_V$    | M  | M  | M  | M  | M  | 10 |
| F$_{VI}$  | M  | M  | M  | M  | 15 | M  |

Next we apply the standard Hungarian assignment algorithm to obtain the initial solution:

Assignment table after first iteration:

|       | F$_I$ | F$_{II}$ | F$_{III}$ | ***F$_{IV}$*** | ***F$_V$*** | ***F$_{VI}$*** |
|-------|-------|----------|-----------|----------|-------|----------|
| ***F$_I$***    | ***M-20*** | ***0***    | ***M-25*** | ***M-20*** | ***M-20*** | ***5***    |
| ***F$_{II}$***  | ***0***    | ***M-10*** | ***0***    | ***M-10*** | ***M-10*** | ***M-10*** |
| F$_{III}$ | M-30 | M-30 | M-35 | ***0***    | M-30 | M-30 |
| F$_{IV}$  | M-40 | M-40 | 5    | ***M-40*** | ***M-40*** | ***0***    |
| F$_V$    | M-10 | M-10 | M-15 | ***M-10*** | ***M-10*** | ***0***    |
| F$_{VI}$  | M-15 | M-15 | M-20 | ***M-15*** | ***0***    | ***M-15*** |

There are two rows and three columns that are covered i.e. *k* = 5. But as this is a 6x6 load matrix, the above solution is sub-optimal. So we make a second iteration:

|       | ***F$_I$*** | F$_{II}$ | ***F$_{III}$*** | ***F$_{IV}$*** | F$_V$ | ***F$_{VI}$*** |
|-------|-------|----------|-----------|----------|-------|----------|
| ***F$_I$***    | ***M-20*** | ***0***    | ***M-25*** | ***M-15*** | ***M-15*** | ***10***   |
| F$_{II}$  | ***0***    | M-10 | ***0***    | ***M-5***  | M-5  | ***M-5***  |
| F$_{III}$ | ***M-35*** | M-35 | ***M-40*** | ***0***    | M-30 | ***M-30*** |
| F$_{IV}$  | ***M-45*** | M-45 | ***0***    | ***M-40*** | M-40 | ***0***    |
| F$_V$    | ***M-15*** | M-15 | ***M-20*** | ***M-10*** | M-10 | ***0***    |
| ***F$_{VI}$***  | ***M-20*** | ***M-20*** | ***M-25*** | ***M-15*** | ***0***    | ***M-15*** |

Now columns $F_I$, $F_{III}$, $F_{IV}$, $F_{VI}$ and rows $F_I$ and $F_{VI}$ are covered i.e. k = 6. As this is a 6x6 load matrix the above solution is optimal.

The optimal assignment table (subject to the 2meters of aisle between adjacent facilities):

|  | $F_I$ | $F_{II}$ | $F_{III}$ | $F_{IV}$ | $F_V$ | $F_{VI}$ |
|---|---|---|---|---|---|---|
| $F_I$ | – | * | – | – | – | – |
| $F_{II}$ | * | – | – | – | – | – |
| $F_{III}$ | – | – | – | * | – | – |
| $F_{IV}$ | – | – | * | – | – | – |
| $F_V$ | – | – | – | – | – | * |
| $F_{VI}$ | – | – | – | – | * | – |

Initial layout of facilities as dictated by the Hungarian assignment algorithm:

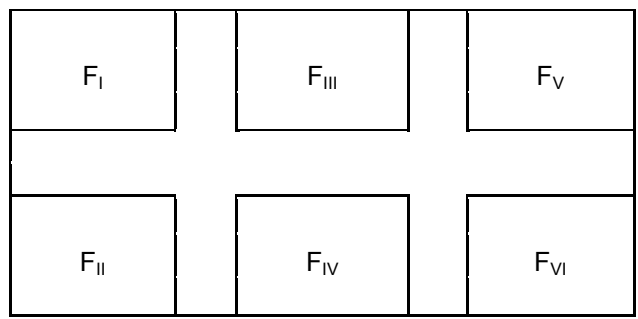

The above layout conforms to the rectangular floor plan of the plant and also places the assigned facilities adjacent to each other with an aisle of 2 meters width between them. Thus $F_I$ is adjacent to $F_{II}$, $F_{III}$ is adjacent to $F_{IV}$ and $F_V$ is adjacent to $F_{VI}$.

Based on the cost information provided in the load-matrix the total cost in terms of load-units for the above layout can be calculated as follows:

L = 2{(20 + 10) + (50 + 30) + (10 + 15)} + (44 x 25) + (22 x 40) + (22 x 15) = 2580.

By feeding the above optimal solution into the CRAFT program the final, the global optimum is found in a single iteration. The final, optimal layout as obtained by CRAFT is as under:

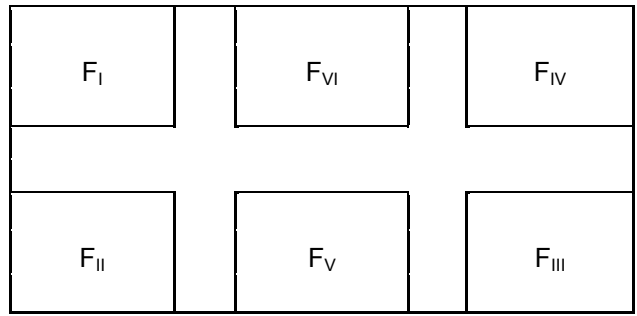

Based on the cost information provided in the load-matrix the total cost in terms of load-units for the optimal layout can be calculated as follows:

L* = 2{(10 + 20) + (15 + 10) + (5 + 30)} + (22 x 25) + (44 x 15) + (22 x 40) = 2360.

Therefore the final solution is an improvement of just 220 load-units over the initial solution! This shows that this initial solution fed into CRAFT is indeed near optimal and can thus ensure a faster convergence.